\documentclass[11pt]{amsart}

\usepackage{amsmath,amsthm,amsfonts,amssymb}
\usepackage[all,cmtip]{xy}
\usepackage{verbatim}
 \usepackage{latexsym}
\newtheorem{thm}{Theorem}
\newtheorem{prop}[thm]{Proposition}
\newtheorem{cor}[thm]{Corollary}
\newtheorem{lemma}[thm]{Lemma}
\newtheorem{definition}[thm]{Definition}

\theoremstyle{definition}
\newtheorem{eg}[thm]{Example}

\theoremstyle{remark}

\newtheorem{case'}{Case}

\numberwithin{equation}{section}
\numberwithin{thm}{section}

\def\({\left(}
\def\){\right)}
\def\[{\left[}
\def\]{\right]}
\def\={\quad = \quad}
\def\+{\quad + \quad}

\def\R{{\mathbb{R}}}
\def\Z{{\mathbb{Z}}}

\def\C{{\mathbb{C}}}

\def \a {\alpha}
\def \b {\beta}

\def \o {\circ}

\def \l {\lambda}

\def \L {\Lambda}

\begin{document}

\title{Gabor Fields and Wavelet Sets for the Heisenberg Group}


\author{Bradley Currey and
Azita Mayeli}

\date{\today}
\thanks{The second author was partially  supported by the Marie Curie Excellence Team Grant MEXT-CT-2004-013477,
 Acronym MAMEBIA, funded by the European Commission.}
\maketitle

 \begin{abstract} We study singly-generated wavelet systems on $\R^2$  that are naturally associated with rank-one wavelet systems on the Heisenberg group $N$. We prove a necessary condition on the generator in order that any such system be a Parseval frame. 
 Given a suitable subset $I$ of the dual of $N$,  we give an explicit construction for Parseval frame wavelets that are associated with $I$.  We say that $g\in L^2(I\times \R)$ is  Gabor field over $I$ if, for a.e. $\l \in I$, $|\l|^{1/2} g(\l,\cdot)$  is the Gabor generator of a Parseval frame  for $L^2(\R)$, and that $I$ is a Heisenberg wavelet set if every Gabor field over $I$ is a Parseval frame (mother-)wavelet for $L^2(\R^2)$. We then show that $I$ is a Heisenberg wavelet  set if and only if $I$ is both translation congruent with a subset of the unit interval and dilation congruent with  the Shannon set.
 
  \end{abstract}

 {\footnotesize
 Keywords and phrases: \textit{Wavelet, Heisenberg group, Gabor frame,
 Parseval frame, multiplicity-free}}

\section{introduction}\label{intro}

The Heisenberg group $N$  is a natural model for the structure of Gabor systems, as well as a domain upon which one can construct wavelet frame systems of functions. There is a substantial body of work that has been concerned with the construction of wavelets on $N$, or more generally, on stratified groups.  In \cite{Lem89}, multiresolution analysis on stratified Lie groups is used to obtain a wavelet orthonormal basis for $L^2(N)$ generated by finitely many $C^N$ wavelets, arising from a generalized spline-surface space. Nearly tight wavelet frames on stratified Lie goups are obtained in \cite{GM06} by a multiplier and sub-Laplacian theory via a single Schwartz  wavelet with many vanishing moments, or a compactly supported smooth wavelet with arbitrarily many vanishing moments. In \cite[Chapter 6]{F}, band-limited subspaces of $L^2(N)$ that admit translation frames are characterized. In \cite{M}, a Shannon-type multiresolution analysis is constructed and used to produce a (single) band-limited  Parseval frame wavelet for $L^2(N)$.  In this paper we are concerned with a class of discrete systems of functions of two variables that involve dilation, as well as translation and modulation, and that arise naturally as wavelet systems of functions on the Heisenberg group. 
 
 Other than in  \cite{F} and  \cite{M}, the group Fourier transform on $L^2(N)$ has played a limited role in the the construction of wavelet frames. This stands in stark contrast with central role played by the Fourier transform in the construction and application of wavelet frames on $\R^n$. In this paper we will see that an explicit application of the group Fourier transform on certain wavelet systems in $L^2(N)$ leads to a natural class of systems in $L^2(\Bbb R^2)$  that are combinations of translations, modulations, and dilations.  
 Our approach is to consider discrete systems whose group Fourier transforms are rank one. The dual $\hat N$ of $N$ is almost-everywhere identified with $\L = \Bbb R\setminus \{0\}$, and hence rank-one  systems can be regarded in a natural way on the Fourier transform side as scalar-valued function systems in $L^2(\L\times \R)$, where the group Plancherel measure $|\l|d\l$ is used on $\L$.  If $\hat \psi = g(\l,t)$ is the generator for such a system, then for each $\l$,  the group translates of $\psi$  become Gabor systems generated by $|\l|^{1/2}g(\l,\cdot)$. Dilation operators that are Fourier transforms of the usual dilations on $L^2(N)$ are added to the mix, and the result is a wavelet system on $L^2(\L\times \R) \simeq L^2(\R^2)$ that is naturally associated with the Heisenberg group. In the present paper, we prove a number of fundamental properties of such systems. We characterize the wavelet systems on $N$ that give rise to wavelet systems in $L^2(\L\times\R)$, and for these systems we prove a necessary condition for frame  generators. We prove a method for construction of wavelet frame systems for $L^2(\L\times \R)$ by considering functions supported on subsets  $I\times \R$  of $\L \times \R$. Given $I \subset \L$, we say that $g\in L^2(I\times \R)$ is a Gabor field over $I$  if, for a.e. $\l \in I$,  the Gabor system generated by  $|\l|^{1/2}g(\l,\cdot)$  is a Parseval frame for $L^2(\R)$. Examples of Gabor fields that generate Parseval frames for $L^2(I\times \R)$ and $L^2(\L\times \R)$ are not difficult to construct, though they can never be orthonormal. We then study conditions under which all Gabor fields over $I$ are also generators for Parseval frames on $L^2(I\times \R)$ and $L^2(\L\times \R)$. The main idea of this paper is that these conditions are solely dependent upon the properties of the set $I$, and that these properties are precisely the analogue of wavelet frame set properties in the Euclidean case. Further study of general wavelet systems if this type is intended for future work.

In Section \ref{group} we introduce invariant multiplicity-free subspaces of $L^2(N)$, as well as multiplicity one subspaces, and we make explicit their natural isomorphism with the function space $L^2(I\times \R)$, where $I$  is the spectrum of the subspace in question. In Section 2 we consider the systems in $L^2(I\times \R)$ that are Fourier transforms of translation systems in multiplicity-free subspaces. Such systems are fields of Gabor systems parametrized by $I$, and in light of basic results of Gabor analysis we observe that a band-limiting condition on $I$  is necessary for the existence of Gabor fields over $I$.  We show  (Propositions \ref{PF on I times R} and \ref{GHS}) that the property that every Gabor field over $I$ generates a Parseval frame for $L^2(I\times \R)$ is equivalent to the property that $I$ is translation congruent with a subset of the unit interval. We then introduce the dilations on $N$ and $L^2(N)$, and show that if $I$ is dilation congruent with the Shannon set, then translation frame generators for $L^2(I\times \R)$ are also wavelet frame generators for $L^2(\L\times\R)$. 

We begin Section \ref{final section} by proving a precise analogue of \cite[Theorem 3.3.1]{D}, showing that if $g$ is the generator of a wavelet frame system for $L^2(\L\times\R)$ (using the integer translation lattice), then $g$ is weakly admissible for the quasi-regular representation of the dilated Heisenberg group (Theorem \ref{nec condition}). With this general result we turn back to band-limited wavelet frame generators. If $I \subset \L$, we say that $g$ is an $I$-wavelet if $g\in L^2(I\times \R)$ and $g$ is the generator of a wavelet frame system for $L^2(\L\times\R)$ using the integer translation lattice. We then show that if every Gabor field $g$ over $I$ is an $I$-wavelet, then $I$ is dilation congruent with the Shannon set, and $g$ is a wavelet frame generator for non-integral translation lattices in $N$ as well. In this context we say that $I$ is a wavelet set if every Gabor field $g$ over $I$ is also an $I$-wavelet. All of the above results allow us to conclude (Theorem \ref{PW sets}) that $I$ is a wavelet  set if and only if $I$ is translation congruent with a subset of the unit interval and $I$ is dilation congruent with the Shannon set.

\section{preliminaries }\label{group}

A system $\mathcal W = \{\psi_j\}_{j\in J}$ of functions in a separable Hilbert space $\mathcal H$ is a \it Parseval frame \rm (or \it normalized tight frame\rm) for $\mathcal H$ if
$$
\|g\|^2 = \sum_{j\in J} \ \left|\langle g, \psi_j\rangle \right|^2
$$
holds for every $g \in \mathcal H$. If $\mathcal W $ is a system in  $\mathcal H$ and $\mathcal W' $ is a system in $\mathcal H'$, we will say that $\mathcal W$ and $\mathcal W'$ are equivalent if there is a unitary isomorphism $V : \mathcal H \rightarrow \mathcal H'$ whose restriction to $\mathcal W$ is a bijection onto $\mathcal W'$.

The term ``set" will always mean ``Lebesgue measurable set" and  ``function'' will mean ``Lebegue measurable function". We shall identify subsets whose symmetric difference has Lebesgue measure zero. We recall the following facts about translation and dilation congruence. First, subsets $I$ and $J$ of $\R$ are translation congruent if there is a bijection $\tau : I\rightarrow J$ and an integer-valued function $k$ on $I$ such that  $\tau(\l)=\l+k(\l)$. $I$ is translation congruent with a subset of the unit interval
$[0,1]$ if and only if the integral translations of $I$ are disjoint. Secondly, $I$ and $J$ are dilation congruent if there is a bijection $\delta: I\rightarrow J$ and an integer valued function $j$ on $I$ such that $\delta(\l) = 2^{j(\l)}\l$.  $I$ is dilation congruent with the Shannon set $I_0 = [-1,-1/2) \cup [1/2,1)$
if and only if $\R = \dot\cup_{j\in \Z}2^{j}I$ as a disjoint union.

We use the following realization of the Heisenberg group $N$: as a topological space $N$ is identified with $\R^3$, and we let $N$ have the group operation
$$
(x_1,x_2,x_3)\cdot (y_1,y_2,y_3) = (x_1 + y_2,x_2+y_2,x_3 + y_3+x_1y_2).
$$
For $a> 0 $ and $x \in N$  define the dilation of $x$ by $a$ with
\begin{align}\label{dilation}
a \cdot x = (a^{1/2}x_1, a^{1/2}x_2, ax_3).
\end{align}
Then $x \mapsto a\cdot x$ is an automorphism of $N$.

Put $\L = \R \setminus \{0\}$. For $x \in N$, $\l \in \L$,  we define the unitary operator $\pi_\l(x)$ on $L^2(\R)$ by
$$
\Bigl(\pi_\l(x)f\Bigr)(t) = e^{2\pi i \l x_3} e^{-2\pi i \l x_2 t} f(t-x_1), \ f \in   L^2(\R).
$$
Then $x \mapsto \pi_\l(x)$ is an irreducible unitary representation of $N$, the so-called Schr\"odinger representation. Recall that the family $\{\pi_\l : \l \in \L\}$ consists of pairwise inequivalent representations and can be regarded as a non-commutative ``frequency domain''  for $N$ in the following sense. For $\phi \in   L^2(N) \cap  L^1(N)$ and $\l \in\L$, the weak operator-valued integral
$$
\mathcal F \phi (\l) = \int_N \phi(x) \pi_\l(x) dx.
$$
defines a trace-class operator on $L^2(\R)$. The Plancherel theorem for $N$ says that $\mathcal F$ extends to a unitary isomorphism of $L^2(N)$ with the Hilbert space $L^2\bigl(\L, \mathcal{HS}(L^2(\R)), |\l|d\l\bigr)$ of all (equivalence classes of) Hilbert-Schmidt operator-valued functions on $\L$. For simplicity we shall use the notation $\mathcal F \phi = \hat\phi$, and we put $\Bbb H = L^2\bigl(\L, \mathcal{HS}(L^2(\R)), |\l|d\l\bigr)$.

 For $x \in N$ define the unitary left translations on $T_x$ on $ L^2(N)$ by $T_x\phi = \phi(x^{-1}\cdot)$; for $\phi \in   L^2(N) \cap  L^1(N)$, we have
\begin{equation}\label{translation to modulation}
\widehat{T_x\phi }(\l) = \pi_\l(x)\hat \phi(\l), \l \in \L.
\end{equation}
Define the unitary dilation operator on $L^2(N)$ by
$$
D_a\phi = \phi(a^{-1} \cdot )\ a^{-1}  , \phi \in   L^2(N).
$$
For any $a> 0$, $\phi\in L^2(N)$, we have (for a.e. $\l$)
\begin{equation}\label{transform of dilation}
\widehat{D_a \phi}(\l) = a\ C_a^{-1} \o \hat \phi(a\l) \o C_a \
\end{equation}
where, for $f \in  L^2(\R)$ and $a>0$,
$$
C_a f(t) = a^{1/4} f(a^{1/2}t).
$$
Subspaces of $L^2(N)$ that are invariant under these operators are described as follows. 

For Hilbert spaces $\mathcal{H}$ and $\mathcal{K}$ we denote the Hilbert space of all Hilbert-Schmidt operators from  $\mathcal{K}$ to  $\mathcal{H}$ by $ \mathcal{H}\otimes \overline{\mathcal{K}}$, and we recall that given bounded linear operators $T_1$ and $T_2$ on $\mathcal{H}$ and $\mathcal{K}$ respectively, then $T_1\otimes T_2$ is the bounded linear operator on $ \mathcal{H}\otimes \overline{\mathcal{K}}$ defined by
\begin{align}
(T_1\otimes T_2)(A) = T_1 A T_2.
\end{align}
If $\{\mathcal H_\a\}$ is a measurable field of Hilbert spaces over
a measurable space $A$ and $\mu$ is a measure on $A$, then the
direct integral $\int^\oplus_A \mathcal H_\a d\mu(\a)$ is the
Hilbert space of measurable vector fields $\{f_\a\}_{\a\in A}$ on
$A$ such that $\int_A \|f_\a\|^2 d\mu(\a) < \infty$. A bounded
linear operator $T$ on $\int^\oplus_A \mathcal H_\a d\mu(\a)$ is
decomposable if for a.e. $\a\in A$, there is a bounded linear
operator $T_\a$ on $\mathcal H_\a$ such that $T\bigl(\{f_\a\}\bigr)
= \{T_\a f_\a\}$, and we write $T = \{T_\a\}_{\a\in A}$. In
particular, we have
$$
\Bbb H = \int^\oplus_\L \mathcal{HS}\bigl(L^2(\R)\bigr)\  |\l|d\l
$$
and for each $x \in N$, $\hat T_x := \mathcal F T_x \mathcal F^{-1} = \{\pi_\l(x)\otimes 1\}_{\l\in\L}$.
 See \cite[Section 7.3, 7.4]{Fo} for further details about the preceding.

Now let $\mathcal H$ be a closed subspace of $L^2(N)$; we say that $\mathcal H$ is  \it translation invariant
 \rm if $T_x(\mathcal H) \subset \mathcal H$ holds for all $x \in N$.  Assume that $\mathcal H$ is translation invariant and let $P$
be the orthogonal projection onto $\mathcal H$. Then $P$ commutes with the operators $T_x, x \in N$, and hence $P$ belongs to
the von-Neumann algebra generated by the right translation operators.  By \cite[Th. 18.8.1]{Di} or \cite[Theorem 3.48]{F}, we have
\begin{equation}\label{proj decomp}
\hat P := \mathcal F P \mathcal F^{-1} = \{1 \otimes \hat P_\l\}_{\l\in\L}
\end{equation}
where, for a.e. $\l$, $\hat P_\l$ is an orthogonal projection on $L^2(\R)$, and the measurable field $\{\hat P_\l\}_{\l\in\L}$
 is unique up to a.e. equality. This means that \begin{equation}\label{inv decomp}
\hat{\mathcal H} := \mathcal F\bigl(\mathcal H\bigr) = \int^\oplus_\L \ L^2(\R) \otimes \overline{\mathcal K}_\l \ |\l| d\l.
\end{equation}
where $\mathcal K_\l =\hat P_\l(L^2(\R))$. Set $m_\mathcal H(\l) = \text{rank}(\hat P_\l), \l \in \L$; the spectrum of $\mathcal H$ is the set $I(\mathcal H) = \text{supp}(m_\mathcal H)$. A closed subspace $\mathcal H$ of $L^2(N)$ is said to be multiplicity-free if it is translation invariant and $m_\mathcal H(\l) \le 1$ a.e.; if  $m_\mathcal H(\l) = 1$ a.e. then we will say that  $\mathcal H$ is \it multiplicity-one \rm.

Let $\mathcal H$ be a multiplicity-free subspace of $L^2(N)$, $I$ the spectrum of $\mathcal H$, $P$ the projection onto $\mathcal H$, and let $\{\hat P_\l\}$ be the associated measurable field of projections as in (\ref{proj decomp}). We have a measurable vector field $e = \{e_\l\}_{\l\in I}$ with each $e_\l$ a unit vector in $L^2(\R)$, such that $\hat P_\l = e_\l \otimes e_\l$ and $\mathcal K_\l = \C e_\l$. Thus (\ref{inv decomp}) can be written as
\begin{equation}\label{mult one int}
\hat{\mathcal H}= \int^\oplus_I \ L^2(\R) \otimes e_\l \ |\l| d\l.
\end{equation}
where $e_\l$ is regarded as an element of $\overline{L^2(\R)}$.
Hence $\mathcal H$ is isomorphic with
\begin{align}\label{reduced-subspace}
\int^\oplus_I \ L^2(\R) \ |\l|d\l
\end{align}
via the isomorphism $V_e$ defined on $\mathcal H$ by
$$
V_e\eta(\l) = \hat\eta(\l)\bigl(e_\l\bigr), \ {\rm a.e.\ } \eta \in \mathcal H.
$$
We identify the direct integral (\ref{reduced-subspace})  with $L^2(I\times\R)$ in the obvious way, where it is understood that $I$ carries the measure $|\l|d\l$. 
Note that the definition of $V_e$ depends upon the choice of the unit vector field $e =\{e_\l\}$, and if we write $\hat\eta= \{ f_\l\otimes e_\l\}$, then  $V_e\eta(\l) = f_\l$. We shall say that $V_e$ is the \it reducing isomorphism \rm for $\mathcal H$ associated with the vector field $\{e_\l\}$. Note that if $e' = \{e'_\l\}$ is another measurable unit vector field for which (\ref{mult one int}) holds, then there is a measurable unitary complex-valued function $c(\l)$ on $I$ such that $e'_\l = c(\l)e_\l$ holds for a.e. $\l$. Hence if $V_e$ and $V_{e'}$ are the associated reducing isomorphisms, then $V'\o V^{-1}$ is just the unitary multiplication operator on (\ref{reduced-subspace})  associated with $c(\l)$. From now on we write $V : \mathcal H \rightarrow L^2(I\times \R)$ when discussing a multiplicity-free subspace with reducing isomorphism $V$.

\vspace{.1in}


\section{Gabor fields, translation systems, and wavelet systems}\label{Gabor}

Let $I \subseteq \L$ and let  $V : \mathcal H \rightarrow L^2(I\times \R)$ be a multiplicity-free subspace of $L^2(N)$ . Let $\Gamma_{\a,\b} = \a \Bbb Z \times \b\Bbb Z \times \Bbb Z$ where $\a$ and $\b$ are positive parameters. We shall refer to the subsets $\Gamma_{\a,\b}$ as lattices and the constants $\a $ and $\b$ as lattice parameters, even though $\Gamma_{\a,\b}$ is not a subgroup of $N$ unless $\a$ and $\b$ are integers. For $(k,l,m) \in \Gamma_{\a,\b}$ we define the unitary operator on $L^2(I\times \R)$ by  $\hat T_{k,l,m} = VT_{(k,l,m)} V^{-1}$, so that 
for $g\in L^2(I\times \R)$, 
$$
\hat T_{k,l,m}g(\l,t) = e^{2\pi i \l m}e^{-2 \pi i \l l t} \ g(\l,t- k).
$$
With $\psi = V^{-1}g$,  the translation system $\mathcal T(\psi, \a,\b) = \{T_{(k,l,m)} \psi : (k,l,m)\in\Gamma_{\a,\b}\}$ is equivalent with the system $
 \widehat{ \mathcal  T}(g, \a,\b) = \{\hat T_{k,l,m}g : (k,l,m)\in\Gamma_{\a,\b}\}.
$
For $\l \in \L$ fixed, $\hat T_{k,l,0}$ defines a unitary operator
on $L^2(\R)$ in the obvious way which we denote by $\hat
T_{k,l}^\l$. For $u \in L^2(\R)$ set   $\mathcal G(u, \a,\b,\l) = \{\hat
T_{k,l}^\l u : k \in \a\Bbb Z, l \in\b\Bbb Z\}$. We will also use the
notation $g_{k,l,m} = \hat T_{k,l,m} g$.

\begin{definition} Let $I \subset \L$ and $g \in L^2(I\times \R)$. We say that $g$ is a Gabor field over $I$  if, for a.e. $\l \in I$, $\mathcal G(|\l|^{1/2} g(\l,\cdot), \a,\b,\l) $ is a Parseval frame  for $L^2(\R)$, for some lattice parameters $\a$ and $\b$. 

\end{definition}

Gabor fields are easily constructed as the following shows. 

\begin{eg} \label{indicator functions} Let $\a$ and $\b$ be any positive numbers, and let  $I$ be a subset of $\L$ that is included in the interval $ [-1/\a\b, 1/\a\b]$.   For  each
 $(\l,t) \in \L \times \R$, set 
 $$
 g(\l,t) =  \b^{1/2} {\bf 1}_{I\times
 [0,\a]}(\l,t).
 $$
 Then $g$ is  a Gabor field over $I$ with lattice parameters $\a$ and $\b$.

\end{eg}

\begin{proof} Fix any $\l \in I $ and for $l \in \b \Bbb Z$ put $e^\l_l(t) = e^{-2\pi i \l l t}$. For $f \in L^2(\R)$, $k\in\a\Bbb Z, l\in\b\Bbb Z$, we have
$$
\begin{aligned}
\langle f, |\l|^{1/2} g_{ k,l,0}(\l,\cdot ) \rangle &= \int_\R \ f(t)\  |\l\b|^{1/2}e^{2\pi i \l l t} {\bf 1}_{[0, \a ]}(t- k) \ dt  \\
&= \left(\int_\R \ f(t+k){\bf 1}_{[0,\a ]}(t) \ |\l\b|^{1/2}e^{2\pi i \l l t} dt\right) e^{2\pi i \l k l}.
\end{aligned}
$$
Since   $\{|\l\b|^{1/2}e^\l_l{|_{[0,1/\b|\l|]}} : l \in \b\Bbb Z\}$ is an orthonormal basis for $L^2([0,1/\b|\l|])$ and $\a \le 1/\b|\l|$,
then  $\{|\l\b|^{1/2}e^\l_l{|_{[0,\a ]}} : l \in \b\Bbb Z\}$ is a Parseval frame for $L^2([0,\a ])$. Hence we have
$$
\begin{aligned}
\|f\|^2 = \sum_{k\in\a\Bbb Z} \| {\bf 1}_{[0,\a ]} \ T_{- k}f\|^2& = \sum_{k\in\a\Bbb Z} \left(\sum_{l\in\b\Bbb Z} \
 \Bigl|\langle {\bf 1}_{[0,\a ]}\ T_{- k}f , |\l\b|^{1/2}e^\l_{l}{|_{[0,\a ]}}\rangle\Bigr|^2 \right)\\
&=\sum_{k\in\a\Bbb Z} \left(\sum_{l\in\b\Bbb Z} \ \Bigl|\langle f ,|\l\b|^{1/2}e^\l_{l}\ T_{ k}{\bf 1}_{[0,\a ]}\rangle
e^{2\pi i\l l k}\Bigr|^2 \right)\\
&= \sum_{k\in\a\Bbb Z, l\in\b\Bbb Z} \ |\langle f, |\l|^{1/2}\ g_{ k,l,0}(\l,\cdot)\rangle |^2.
\end{aligned}
$$
\end{proof}

Now suppose that  $g$ is the function of the preceding example and also that  $I$ is translation congruent with a subset of the unit interval. Then the system $\hat{\mathcal T}(g,\a,\b)$ is a Parseval frame for $L^2(I\times \R)$, as the following result shows.

\begin{prop} \label{PF on I times R} Suppose that  $I$ is translation congruent with a subset of the unit interval, and fix lattice parameters $\a$ and $\b$. For each $g \in L^2(I\times \R)$,  the conditions

\vspace{.1in}
\noindent
(i) $g$ is a Gabor field over $I$ with lattice parameters $\a$ and $\b$, and 

\vspace{.1in}
\noindent
(ii)  $\widehat{ \mathcal  T}(g, \a,\b)$  is a Parseval frame for $L^2(I\times \R)$,

\vspace{.1in}
\noindent
are equivalent. 

\end{prop}

\begin{proof}Let $g \in L^2(I\times \R)$ and suppose that $g$ defines a Gabor field over $I$ with lattice parameters $\a$ and $\b$. For    $f  \in L^2(I \times \R)$ we have $f(\l,\cdot)$ is square integrable on $\R$ for a.e.$\l$, so
$$
\|f(\l,\cdot)\|_{L^2(\R)}^2 = \sum_{k\in\a\Bbb Z, l\in \b\Bbb Z} \ |F_{k,l}(\l)|^2.
$$
holds for a.e. $\l \in I$, where
$$
F_{k,l}(\l) =  \langle f(\l,\cdot),|\l|^{1/2}g_{k,l,0}(\l,\cdot)\rangle.
$$
Hence
\begin{equation}\label{f norm}
\|f\|^2 = \int_I \ \|f(\l,\cdot)\|_{L^2(\R)}^2 \ |\l|d\l = \sum_{k\in\a\Bbb Z, l\in \b\Bbb Z} \   \int_I \ |F_{k,l}(\l)|^2 \ |\l| d\l.
\end{equation}
and in particular, each $F_{k,l}$ belongs to $L^2(I,|\l|d\l)$.
Now since $\{e^{2\pi i\l m}|_{[0,1]}: m \in \Bbb Z\}$ is an orthonormal basis for $L^2([0,1],d\l)$ and $I$ is  translation
congruent to a measurable subset of $[0,1]$, then
 $\{e^{2\pi i\l m}{|_I} : m \in \Bbb Z\}$ is a  Parseval frame for $L^2(I,d\l)$. Hence
we have
$$
\begin{aligned}
\|F_{k,l}\|_{L^2(I)}^2 &= \int_I \ \left| \langle f(\l,\cdot), |\l|^{1/2} g_{k,l,0}(\l,\cdot) \rangle\right|^2 \ |\l|d\l \\
&= \int_I \ \left|\langle f(\l,\cdot), g_{k,l,0}(\l,\cdot) \rangle |\l| \ \right|^2 \ d\l \\
&= \sum_{m \in \Bbb Z} \ \left|\int_I \  \langle f(\l,\cdot), g_{k,l,0}(\l,\cdot) \rangle |\l| \ e^{-2\pi i \l m} \  d\l \right|^2 \\
&= \sum_{m \in \Bbb Z} \ \left| \int_I \ \langle f(\l,\cdot), g_{k,l,m}(\l,\cdot) \rangle |\l| \ d\l \right|^2 \\
&=  \sum_{m \in \Bbb Z} \ \left|\int_{I\times \R} \ f(\l,t) \overline{g_{k,l,m}(\l,t)} \ |\l| d\l dt\right|^2.
\end{aligned}
$$
Combining the preceding with (\ref{f norm}) we get
$$
\|f\|^2 =  \sum_{k\in\a\Bbb Z, l\in \b\Bbb Z,m\in\Bbb Z} \ \left| \int_{I\times \R} \ f(\l,t) \overline{g_{k,l,m}(\l,t)} \ |\l| d\l dt\right|^2.
$$

On the other hand, suppose that $\widehat{ \mathcal  T}(g, \a,\b)$  is a Parseval frame for $L^2(I\times \R)$. For any $f \in L^2(I \times \R)$, we have
\begin{equation}\label{f norm 2}
\|f\|^2 = \sum_{k\in\a\Bbb Z, l\in \b\Bbb Z} \ \sum_m \left| \int_I \langle f(\l,\cdot),g_{k,l,m}(\l,\cdot)\rangle |\l|d\l\right|^2.
\end{equation}
Since $I$ is translation congruent to a subset of $[0,1]$, functions on $I$ that are square-integrable with respect to Lebesgue measure have a Fourier series expansion in terms of the exponentials $e^{2 \pi i m \l}$. Now (\ref{f norm 2}) shows that $\l \mapsto \langle f(\l,\cdot),g_{k,l,0}(\l,\cdot)\rangle |\l|$ is integrable with square integrable Fourier coefficients
$$
\int_I \langle f(\l,\cdot),g_{k,l,0}(\l,\cdot)\rangle |\l| \ e^{-2\pi i \l m} d\l  = \int_I \langle f(\l,\cdot),g_{k,l,m}(\l,\cdot)\rangle |\l|d\l.
$$
Hence $\l \mapsto \langle f(\l,\cdot),g_{k,l,0}(\l,\cdot)\rangle |\l| $ belongs to $L^2(I,d\l)$ and $F_{k,l}(\l) := \langle f(\l,\cdot),|\l|^{1/2}g_{k,l,0}(\l,\cdot)\rangle$ belongs to $L^2(I,|\l|d\l)$. Moreover (\ref{f norm 2}) and the Parseval identity for Fourier series shows that
\begin{equation} \label{f norm 3}
\|f\|^2 = \sum_{k\in\a\Bbb Z, l\in \b\Bbb Z} \ \|F_{k,l}\|^2.
\end{equation}
Now we claim that for $f \in L^2(I\times \R)$, there is a conull subset $C(f)$ of $I$ such that for all $\l \in C(f)$,
$$
\int_\R \ |f(\l,t)|^2 \ dt = \sum_{k\in\a\Bbb Z, l\in \b\Bbb Z} \ |F_{k,l}(\l)|^2
$$
holds. Let $B$ be any Lebesgue measurable subset of $I$, and let $f_B = f \bold 1_{B \times \R}$. Then
$$
\langle f_B(\l,\cdot),|\l|^{1/2}g_{k,l,0}(\l,\cdot)\rangle = F_{k,l}(\l)\bold 1_B(\l)
$$
and so by (\ref{f norm 3}),
\begin{align}\notag
\int_B \int_\R \ |f(\l,t)|^2 \ dt\  |\l|d\l &= \int_I \int_\R \ |f_B(\l,t)|^2 \ dt \ |\l|d\l \\\notag
&= \|f_B\|^2 = \int_B \ \sum_{k\in\a\Bbb Z, l\in \b\Bbb Z} \ |F_{k,l}(\l)|^2 \ |\l|d\l.
\end{align}
The claim follows. Now let $\mathcal D$ be a countable subset of $L^2(I\times \R)$ such that for each
$\l \in I$, $\{f(\l,\cdot) : f \in \mathcal D\}$ is dense in $L^2(\R)$. (For instance, one can take a
 countable dense subset of $L^2(\R)$, multiplied by a fixed nowhere vanishing function in $L^2(\L)$.)
 Let $C = \cap_{f \in \mathcal D} \ C(f)$ where $C(f)$ is as in the preceding claim. For each $\l \in C$ define $T_\l : L^2(\R) \rightarrow l^2(\a\Bbb Z\times\b\Bbb Z)$ by
$$
T_\l(h) = \left( \langle h,|\l|^{1/2}g_{k,l,0}(\l,\cdot)\rangle \right)_{k,l}.
$$
For $f \in \mathcal D$ we have
\begin{align}\notag
\|T_\l(f(\l,\cdot))\|^2 &= \sum_{k\in\a\Bbb Z, l\in \b\Bbb Z} \ \left| \langle f(\l,\cdot),|\l|^{1/2}g_{k,l,0}
(\l,\cdot)\rangle\right|^2 \\\notag
&= \sum_{k\in\a\Bbb Z, l\in \b\Bbb Z}\ |F_{k,l}(\l)|^2 = \|f(\l,\cdot)\|^2
\end{align}
and since $\{f(\l,\cdot) : f \in \mathcal D\}$ is dense in $L^2(\R)$, it follows that $T_\l$ is an isometry.
This means that for each $\l \in C$, $\{ |\l|^{1/2} g_{k,l,0}(\l,\cdot) : k,l \in \Bbb Z\}$ is a Parseval frame for $L^2(\R)$, and the proof is finished.

\end{proof}

The assumption that $I$ is included in $[-1/\a\b, 1/\a\b]$ in Example \ref{indicator functions} is necessary: from \cite[Corollary 7.3.2, Corollary 7.5.1]{G}, we see that if $g$ defines a Gabor field over $I$ with lattice parameters $\a$ and $\b$, then 
\begin{equation}\label{bandlimit}
\||\l|^{1/2}g(\l,\cdot)\|^2 = |\l|\a\b \le 1
\end{equation}
holds for a.e. $\l \in I$. From now on we shall assume then that $I \subset [-1/\a\b, 1/\a\b]$. We first consider the case where $\a = \b = 1$.

\begin{definition} We say that a subset $I$ of $[-1,1]$ is a Gabor-Heisenberg (GH) frame set  if  the following holds.  For every $g \in L^2(I\times \R)$, if $g$ is a Gabor field over $I$ with lattice parameters $\a = \b = 1$, then the Heisenberg system $\widehat{ \mathcal  T}(g, 1,1)$ is a Parseval frame for $L^2(I\times \R)$. 
\end{definition}

We remark that the preceding notion is analogous to the conditions of \cite[Lemma 5.3]{HL} for the situation when the domain is $\R$. In that case $\L = \hat \R$ -- the frequency domain for the Euclidean Fourier transform -- and integer translations on the group $\R$ become modulation on the Fourier transform side, rather than Gabor systems. There one could say that $I$ is a ``modulation frame set'' if it satisfies the following. If $g \in L^2(I)$ and the constant $g(\l)$ is a Parseval frame for $\C$ for a.e. $\l \in I$ (meaning only that $|g(\l)| = 1$ a.e. on $I$), then $e^{2\pi i \l m}g(\l)$ is a Parseval frame for $L^2(I)$. A consequence of \cite[Lemma 5.3]{HL}  in that case (allowing for differences in scaling) is that $I$ is a modulation frame set if and only if $I$ is translation congruent with a subset of the unit interval. We show that an analogous result holds here.

\begin{prop} \label{GHS} Let $I$ be a measurable subset of $[-1,1]$. Then the following are equivalent. 

\vspace{.1in}
\noindent
(i) $I$ is a GH set.

\vspace{.1in}
\noindent
(ii) $I$ is translation congruent with a subset of the unit interval.

\end{prop}

\begin{proof} In light of Proposition \ref{PF on I times R}, we need only show that (i) implies (ii). Suppose that $I$ is not translation congruent with a subset of the unit interval. Since $I$ is incuded in $[-1,1]$, then we have a measurable non-null set $E \subset I \cap [0,1]$ such that $E - 1 \subset I \cap [-1,0]$. Define $g(\l,t) = \bold 1_{[-1,0]}(t)$ for $\l \in I \cap [0,1]$ and $g(\l,t) = \bold 1_{[0,1]}(t)$ for $\l \in I \cap [-1,0]$. It is clear that $g$ defines a Gabor field over $I$. 
Set $\eta(\l,t) = \text{sign}(\l) \bigl( \eta_1(\l,t) + \eta_2(\l,t)\bigr)$ where
$$
\eta_1(\l,t) = \bold 1_E(\l) \bold 1_{(\l-1)E}(t), \ \ \eta_2(\l,t) = \bold 1_{E-1}(\l) \bold 1_{(\l+1)E}(t). 
$$
We claim that $\langle \eta, g_{k,l,m}\rangle = 0$ for all $k, l ,m$. 

We have 
$$
\begin{aligned}
\langle  \eta, g_{k,l,m}\rangle &= \int\int \eta_1(\l,t) e^{-2\pi i \l m} e^{2\pi i l \l t} \l dt d\l +  \int\int \eta_2(\l,t) e^{-2\pi i \l m} e^{2\pi i l \l t} \l dt d\l \\
&=  \int \bold 1_E(\l) e^{-2\pi i \l m} \int \bold 1_{(\l-1)E}(t) e^{2\pi i l \l t} \bold 1_{[-1,0]}(t - k) dt \l d\l  \ \ +\\
& \ \ \ \ \ \ \ \ \ \ \ \ \ \int \bold 1_{E-1}(\l) e^{-2\pi i \l m} \int \bold 1_{(\l+1)E}(t) e^{2\pi i l \l t} \bold 1_{[0,1]}(t - k) dt \l d\l.
\end{aligned}
$$
Now for $\l \in E$, $0 \le \l \le 1$ so $-1 \le \l - 1 \le 0$ and $(\l-1)E \subset [-1,0]$. Hence $\bold 1_{(\l-1)E}(t) \bold 1_{[-1,0]}(t-k) = 0$ unless $k = 0$. Similarly,  for $\l \in E-1$, $\bold 1_{(\l+1)E}(t) \bold 1_{[0,1]}(t-k) = 0$ unless $k = 0$. Hence $\langle  \eta, g_{k,l,m}\rangle = 0$ unless $k = 0$. Now when $k = 0$ we have
$$
\begin{aligned}
\langle  \eta, g_{0,l,m}\rangle &=  \int \bold 1_E(\l) e^{-2\pi i \l m} \ \l  \int \bold 1_{(\l-1)E}(t) e^{2\pi i l \l t} dt d\l \ \ + \\
& \ \ \ \ \ \ \ \ \  \int \bold 1_{E-1}(\l) e^{-2\pi i \l m}  \ \l  \int \bold 1_{(\l+1)E}(t) e^{2\pi i l \l t} dt d\l
\end{aligned}
 $$
 Change $\l \mapsto \l - 1$ in the second term of the preceding and we get
 $$
 \langle  \eta, g_{0,l,m}\rangle =  \int \bold 1_E(\l) e^{-2\pi i \l m}\left(\int F_1(\l,t) dt + \int F_2(\l,t) dt \right) d\l
 $$
where $F_1(\l,t) = \l\ \bold 1_{(\l-1)E}(t) e^{2\pi i l\l t}$ and $F_2(\l,t) = (\l-1)\ \bold 1_{\l E}(t) e^{2\pi i l(\l-1)t}$. 
But for each $\l \in E$ (except $0$), $\l > 0$ while $\l-1 < 0$, so
$$
\int F_1(\l,t) dt \overset{(s = \l t)}= \int \bold 1_{(\l-1)E}\bold 1(s/\l) e^{2\pi i l s} ds = \int \bold 1_{\l(\l-1)E}(s)e^{2\pi i l s} ds 
$$
while
$$
\int F_2(\l,t) dt \overset{(s = (\l-1)t)}= -\int \bold 1_{\l E}\bold 1(s/(\l-1)) e^{2\pi i l s} ds = - \int \bold 1_{\l(\l-1)E}(s)e^{2\pi i l s} ds
$$
so that $\langle  \eta, g_{0,l,m}\rangle = 0 $ also.

\end{proof}

We state the preceding in terms of multiplicity-free subspaces also.

\begin{cor} Let $V : \mathcal H \rightarrow L^2(I\times \R)$ be a multiplicity-free subspace of $L^2(N)$ with $I \subset [-1,1]$. Then the following are equivalent.

\vspace{.1in}
\noindent
(i) For every Gabor field $g$ over $I$, $\mathcal T(V^{-1}g, 1,1)$ is a Parseval frame for $\mathcal H$. 

\vspace{.1in}
\noindent
(ii) $I$ is translation congruent with a subset of the unit interval.
\end{cor}

Next we consider the dilation operators on the subspaces $L^2(I\times \R)$ that arise from dilations $D_a$ via reducing isomophisms. Let $\mathcal H$ be a multiplicity-one subspace of $L^2(N)$. If $D_a(\mathcal H) \subset \mathcal H$ holds for all $a > 0$ then we will say that $\mathcal H$ is  an invariant multiplicity-one subspace of $L^2(N)$. In this case the reducing isomorphisms are semi-invariant for the operators $C_a, a > 0$. 

\begin{lemma} Let $\mathcal H$ be an invariant multiplicity-one subspace of $L^2(N)$, and choose a measurable unit vector field $\{e_\l\}_{\l\in\L}$ such that (\ref{mult one int}) holds. Then for each $a>0$, there is associated with $\mathcal H$ and $\{e_\l\}$ a measurable unitary function $\chi_a: \L \rightarrow \Bbb T$ such that
\begin{equation}\label{chi func}
C_a e_\l = \chi_a(\l) e_{a\l}
\end{equation}
holds for a.e. $\l$. Moreover, the map $a\mapsto \chi_a$ is a homomorphism.

\end{lemma}

\begin{proof} Fix $a> 0$. For each $\l$ put $\hat P_\l = e_\l \otimes e_\l$, so that (\ref{proj decomp}) holds. Set $\hat Q_\l := C_a^{-1} \hat P_{a\l} C_a, \l \in\L$, and define the operator $Q$ on $L^2(N)$ by
$$
\hat Q := \mathcal F Q \mathcal F^{-1} = \{1 \otimes \hat Q_\l \}_{\l\in\L}. 
$$
Using the relation (\ref{transform of dilation}),  it is straightforward to check that $Q\eta = \eta$ for all $\eta\in\mathcal H$ and $Q\eta = 0$ for all $\eta\in\mathcal H^\perp$; hence $Q = P$ and by the uniqueness of the decomposition (\ref{proj decomp}), we have a co-null subset $E_a$ of $\L$ such that for all $\l \in E_a$, $\hat Q_\l = \hat P_\l $. Let $S = \{ \l \in\L :\hat P_\l \ne 0\}$ and fix $\l \in E_a\cap S$. Then for each $f \in L^2(\R)$ we have
$$
\langle f, C_a^{-1} e_{a\l}\rangle C_a^{-1} e_{a\l} = \hat Q_\l f = \hat P_\l f = \langle f, e_{\l}\rangle  e_{\l}.
$$
Set $\chi_a(\l) = \langle C_ae_\l, e_{a\l}\rangle$ for $\l \in E_a\cap S$ and $\chi_a(\l) = 0$ otherwise. Putting $f = e_\l$ in the preceding we find that for $\l \in E_a\cap S$,
$$
C_a e_\l = \chi_a(\l) e_{a\l}.
$$
Since $e_\l$ is a measurable vector field and $\chi_a(\l) = \langle C_ae_\l,e_{a\l}\rangle$ holds for a.e. $\l$, then (by putting $\chi_a(\l) = 0$ for all $\l \notin E_a\cap S$), $\chi_a$ determines a unique (up to a.e. equality) measurable function. Now let $a $ and $b$ be positive numbers. For all $\l$ belonging to the co-null set $E_{ab} \cap E_b \cap b^{-1}E_a \cap S$, we have
$$
\chi_{ab}(\l)e_{ab\l} = C_{ab}e_\l = C_a C_b e_\l = C_a \chi_b(\l) e_{b\l} = \chi_a(\l)\chi_b(\l) e_{ab\l}.
$$

\end{proof}


For $j \in \Bbb Z$, set $D_j = D_{2^j}$. Given $\psi \in L^2(N)$, we consider the wavelet system
\begin{equation}\label{wavelet}
\mathcal U(\psi,\a,\b) =  \{D_jT_{(k,l,m)}\psi  :j\in\Bbb Z, k\in\a\Bbb Z, l\in\b\Bbb Z, m\in\Bbb Z\}.
\end{equation}
Suppose we are given a family $\chi = \{\chi_a : a > 0\}$ of measurable unitary functions on $\L$ such that $\chi_{a^{-1}}
= \overline{\chi_a}$. Define the unitary dilation operators $\hat D^\chi_a, a > 0$ on $L^2(\L\times \R)$ by $\hat D^\chi_a f(\l,t) =
\chi_a(\l) f(a\l,a^{-1/2}t)\ a^{3/4} $, and put $\hat D^\chi_j = \hat D^\chi_{2^j}, j \in \Bbb Z$. Given $g \in L^2(\L\times \R)$ and lattice parameters $\a$ and $\b$, we put 
$$
\widehat{\mathcal U}(g, \a,\b, \chi) = \{\hat D^\chi_j \hat T_{k,l,m}g : j\in\Bbb Z, k\in\a\Bbb Z, l\in\b\Bbb Z, m\in\Bbb Z\}.
$$
The following
shows that the preceding system is equivalent with the wavelet system (\ref{wavelet}).

\begin{prop} \label{translation dilation relation} Let $\mathcal H$ be an invariant multiplicity-one subspace of $L^2(N)$, and let
 $V : \mathcal H \rightarrow L^2(\L\times \R)$ be a reducing isomorphism with associated multiplier $\chi$. Then for each $a > 0$ and $x\in N$
 the following diagram commutes:
\begin{equation}\label{ComDiag}
 \begin{aligned}
    {
  \xymatrix{
\mathcal H \ar[d]_V \ar[r]^{T_x} &\mathcal H   \ar[r]^{D_a}  &\mathcal H   \ar[d]^V\\
L^2(\L\times \R)  \ar[r]^{\hat T_x} &L^2(\L\times \R) \ar[r]^{\hat D_a^\chi}  & L^2(\L\times \R)  }}
 \end{aligned}
 \end{equation}
 and hence the system $\mathcal U(\psi, \a,\b)$ is equivalent to  $\widehat{\mathcal U}(V\psi, \a,\b,\chi)$
   via the isomorphism $V$.
 \end{prop}

\begin{proof}   We need only show how the isomorphism $V$ interwines
with dilation and translation operators.  For this,
let $\eta \in \mathcal H$. For each $x  \in N$ we have
$$
\bigl(VT_x\eta\bigr)(\l) = \widehat{T_x\eta}(\l)(e_\l) = \pi_\l(x) \bigl( \hat\eta(\l)(e_\l)\bigr) = \hat T_xV\eta(\l).
$$
For $a> 0$, we again use equation (\ref{transform of dilation}); for a.e. $\l$,
\begin{align}\notag
\bigl(VD_a\eta\bigr)(\l) &= \widehat{D_a\eta}(\l)(e_\l) = C_a^{-1}\hat\eta(a\l) C_a e_\l \ a \\\notag
&= C_a^{-1} \hat\eta(a\l)\bigl(\chi_a(\l) e_{a\l}\bigr) \ a \\\notag
&= \chi_a(\l) C_a^{-1}\bigl(V\eta\bigr)(a\l) \ a \\\notag
&= \bigl( \hat D_a^\chi V\eta\bigr)(\l).
\end{align}
\end{proof}

Now suppose that  $I$ is a subset of $\L$ that meets each dilation orbit $\{2^j\l : j \in \Bbb Z\}$ at exactly one point. Recall that this means that  $I$  is dilation congruent with the Shannon set, and that $\R = \dot\cup_{j\in \Z}2^jI$. For $f \in L^2(\L\times \R)$ we have
$$
\begin{aligned}
\int_{\L \times \R} \ |f(\l,t)|^2  |\l|d\l dt &= \int_\L \ \sum_j {\bf 1}_{2^{-j}I }(\l)\ \int_\R \ \left|  f(\l,t)\right|^2   dt |\l|d\l \\
&= \int_{\L \times \R} \  \sum_j {\bf 1}_{2^{-j}I }(\l)\ \left| f(\l,2^{j/2}t)2^{j/4}\right|^2 \   dt |\l|d\l\\
&= \sum_{j\in\Bbb Z} \ \int_{\L \times \R} \ \left| {\bf 1}_{I \times \R}(\l,t) \ \chi_{-j}(\l) \ f(2^{-j}\l,2^{j/2}t)\  2^{-3j/4}\right|^2 \
  dt |\l|d\l,
\end{aligned}
$$
that is,
\begin{equation}\label{dilation integral}
\|f\|^2 = \sum_{j \in \Bbb Z} \ \|{\bf 1}_{I \times \R} \ \hat D^\chi_{-j}f\|^2
\end{equation}
In light of the preceding as well as Proposition \ref{translation dilation relation} the following is immediate.

\begin{prop}\label{the concrete frame}  Let $\mathcal H$ be any invariant multiplicity-one subspace of $L^2(N)$ with reducing isomorphism $V : \mathcal H \rightarrow L^2(\L\times \R)$ and  let $\a$ and $\b$ be any lattice parameters.  Suppose that $I$ is dilation congruent with the Shannon set and let $g\in L^2(I\times \R)$, $\psi = V^{-1}g$. If  $\hat{\mathcal T}(g,\a,\b)$ is a Parseval frame for $L^2(I\times\R)$, then the wavelet systems  $\mathcal U(\psi, \a,\b)$ and 
$ \widehat{\mathcal U}(g, \a, \b, \chi)$ are a Parseval frames for $\mathcal H$ and $L^2(\L\times\R)$, respectively. 
\end{prop}

\begin{proof} Let $f \in  L^2(\L\times\R)$; since $g$ has support in $I  \times \R$ we have
$$
\| {\bf 1}_{I \times \R} \ \hat D^\chi_{-j}f \|^2= \sum_{k,l,m} \left| \langle \hat D^\chi_{-j}f, g_{k, l,m}\rangle \right|^2 = \sum_{ k,l,m} \left| \langle f, \hat D^\chi_{j} g_{ k, l,m}\rangle \right|^2.
$$
Now combine the preceding with equation (\ref{dilation integral}).

\end{proof}

We apply the preceding to Example \ref{indicator functions}.

 \begin{cor} \label{indicator frames} Suppose that $I$ is both translation congruent with a subset of the unit interval and is dilation congruent with the Shannon set. Set $g =\b^{1/2}\bold 1_{I\times[0,\a]}$ and put $\psi = V^{-1}g$.  Then the systems  $\mathcal U(\psi, \a,\b)$ and 
$ \widehat{\mathcal U}(g, \a, \b, \chi)$ are a Parseval frames for $\mathcal H$ and $L^2(\L\times\R)$, respectively. 
 \end{cor}

 \section{A necessary condition for wavelet generators and wavelet sets}\label{final section}

 Let $G$ denote the semi-direct product $N \rtimes H$ of $N$ by the multiplicative group $H$ of positive real numbers, where $H$ acts on $N$ by the dilations defined in (\ref{dilation}). Recall that a left Haar measure is given on $G$ by $a^{-3} dx da$.
Let $\tau$ be the quasiregular representation of $G$ acting in $L^2(N)$: for $\phi \in L^2(N)$,
$$
\tau(x,a)\phi = T_xD_a\phi =  \phi(a^{-1}x^{-1}\cdot) a^{-1}, x \in N, a \in H.
$$
Observe that the invariant multiplicity one subspaces of $L^2(N)$ are $\tau$-invariant; indeed, the restriction of $\tau$ to $\mathcal H$ is a direct sum of two (non-isomorphic) irreducible representations $\rho_+$ and $\rho_-$ of $G$.

Now let $\mathcal H$ be an invariant multiplicity one subspace and choose a reducing isomorphism $V$ for $\mathcal H$ with multiplier $\chi$. Then we have seen that
$V\tau(x,a)|_{\mathcal H} V^{-1} = \hat T_x \hat D_a^\chi.$ (See diagram (\ref{ComDiag}).) 
Now as is well-known (see for example \cite{C3}), there are non-zero functions $\psi\in L^2(N)$ such that for every $\phi \in L^2(N)$, $(x,a) \mapsto \langle \phi,\tau(x,a)\psi\rangle$ belongs to $L^2(G)$. It is easily seen that projections of such functions to $\mathcal H$ have the same property. Via the isomorphism $V$, it follows that we can choose a non-zero  $h \in L^2(\L\times \R)$ such that for every $f \in L^2(\L\times \R)$,
$$
\int_G \ \left| \langle f, \hat T_x \hat D_a^\chi h\rangle \right|^2 a^{-3} dxda < \infty,
$$
and the linear function $W_h : L^2(\L\times \R) \rightarrow L^2(G)$ defined by $W_h(f)(x,a) =  \langle f, \hat T_x \hat D_a^\chi h\rangle$ is bounded. Now it is straightforward to construct an invariant multiplicity one subspace with reducing isomorphism whose multiplier $\chi$ is trivial: $\chi_a(\l) = 1$ for all $a> 0, \l \in \L$ (see for example \cite{C3} or \cite{M}). For the purposes of this section we assume that this is the case, and we put $\hat D_a = \hat D_a^\chi$.

Define the linear operator $C$ on $L^2(\L\times \R)$ by the weak operator-valued integral
$$
C(f) = \int_{a = 1}^2 \ \int_N \ W_h(f)(x,a) \hat T_x \hat D_a h \ c(x,a) \ a^{-3} dxda
$$
where $c(x,a)$ is a positive, integrable function on $G$ with supp$(c) \subset N \times [1,2]$.

\begin{lemma} \label{trace class} $C$ is a positive trace-class operator and
$$
{\rm Tr\ }C = \|h\|^2 \int_1^2 \int_N c(x,a) \ a^{-3} \ dxda.
$$
\end{lemma}

\begin{proof} For $f\in L^2(\L\times \R)$ and $(x,a)\in G$ we have
\begin{align}\notag
 \langle Cf,f\rangle &= \int_1^2 \int_N \ W_h(f)(x,a) \langle \hat T_x \hat D_a h,f\rangle c(x,a) a^{-3}dxda\\\notag
 & = \int \int \ |W_h(f)(x,a)|^2 c(x,a) a^{-3}dxda
\end{align}
showing that $C$ is positive. Let $\{f_n\}$ be an orthonormal basis for $L^2(\L\times \R)$. Then for each $(x,a)\in G$ we have
$$
\sum_n \ |W_h(f_n)(x,a)|^2 = \sum_n  \ |\langle f_n, \hat T_x \hat D_a h\rangle |^2 = \|\hat T_x \hat D_a h\|^2 = \|h\|^2
$$
and so
$$
\begin{aligned}
\sum_n \ \langle Cf_n,f_n\rangle &= \int_1^2 \int_N  \left(\sum_n \ |W_h(f_n)(x,a)|^2\right) c(x,a) a^{-3}dxda \\
&=  \|h\|^2 \ \int_1^2\int_N c(x,a) \ a^{-3} dxda.
\end{aligned}
$$

\end{proof}

We are now ready to prove the following result for arbitrary wavelet frame generators in $L^2(\L\times \R)$, which is an analogue of \cite[Theorem 3.3.1]{D}.

\begin{thm}\label{nec condition} Suppose that  \ $\widehat{\mathcal U}(g, \a,\b,1)$ is a frame for $L^2(\L\times \R)$ with frame bounds $A,B$. Then
$$
A \ \a\b\ln 2  \le \int_0^\infty  \int_\R \ |g(\l,t)|^2 \ \l^{-1} dt d\l \le  B \ \a\b\ln 2
$$
 and
 $$
A \ \a\b\ln 2 \le \int_{-\infty}^0 \int_\R \ |g(\l,t)|^2\   |\l|^{-1}dt d\l  \le B \ \a\b\ln 2.
$$
\end{thm}

\begin{proof} We prove the first equality; the proof of the second is similar. Choose a function $h$ as above so that $W_h$ is bounded, and so that supp$(h) \subset (0,\infty) \times \R$.
 For ease of notation we write $g_{jklm} :=\hat D_j \hat T_{k,l,m}g$. We compute that
$$
\begin{aligned}
\langle C g_{jklm},g_{jklm}\rangle &= \int_1^2 \int_N \ \left|\langle g_{jklm},\hat T_x\hat D_a h\rangle \right|^2 \ c(x,a) \ a^{-3} \ dx da \\
&= \int_1^2 \int_N \ \left|\langle g,\hat T_{k,l,m}^{-1} \hat T_{2^{-j}\cdot x}\hat D_{2^{-j}a} h\rangle \right|^2 \ c(x,a) \ a^{-3} \ dx da\\
&=  \int_1^2 \int_N \ \left|\langle g,\hat T_{(k,l,m)^{-1}x}\hat D_{2^{-j}a} h\rangle \right|^2 \ c((2^j\cdot x),a) \ a^{-3} 2^{2j} \ dx da\\
&= \int_{2^{-j}}^{2^{-j+1}} \int_N \ \left|\langle g, \hat T_x \hat D_{a} h\rangle \right|^2 \ c(2^j\cdot ((k,l,m)x),2^ja) \ a^{-3} \ dx da.
\end{aligned}
$$

Now let $\delta$ be any positive number, and let $c(x,a) $ be defined for $1 \le a \le 2$ by
$$
c(x,a) = w\left(\frac{|x_1|}{\sqrt{a}}\right)w\left(\frac{|x_2|}{\sqrt{a}}\right)w\left(\frac{|x_3|}{a}\right)
$$
where $w(s) = w_\delta(s) = \delta \ e^{-\pi(\delta  s)^2}$, and $c = 0$ otherwise. For $2^{-j} \le a \le 2^{-j+1}$,
$$
\begin{aligned}
c(2^j\cdot &((k,l,m)x),2^ja) \\ &=w\left(\frac{|2^{j/2}(x_1+k)|}{\sqrt{2^{j}a}}\right)w\left(\frac{|2^{j/2}(x_2+l)|}{\sqrt{2^{j}a}}\right)w\left(\frac{|2^{j}(x_3+m+kx_2)|}{2^{j}a}\right)\\
&= w\left(\frac{|x_1+k|}{\sqrt{a}}\right)w\left(\frac{|x_2+l|}{\sqrt{a}}\right)w\left(\frac{|x_3+m+kx_2)|}{a}\right).
\end{aligned}
$$
So we can write
$$
\begin{aligned}
T_\delta&:= \sum_{j,k,l,m} \langle C g_{jklm},g_{jklm}\rangle\\
&=  \int_{0}^{\infty} \int_N \ \left|\langle g, \hat T_{x} \hat D_{a} h\rangle \right|^2 \ S_\delta(x,a)  \ a^{-3} dx da
\end{aligned}
$$
where
$$
S_\delta(x,a) =  \sum_{k,l,m} \  w\left(\frac{|x_1+k|}{\sqrt{a}}\right)w\left(\frac{|x_2+l|}{\sqrt{a}}\right)w\left(\frac{|x_3+m+kx_2|}{a}\right).
$$
 Now since $C$ is a positive trace class operator on $L^2(\L\times \R)$ and $\{g_{jklm} : j\in \Bbb Z, k\in\a\Bbb Z, l\in\b\Bbb Z, m\in\Bbb Z\}$ is a frame for $L^2(\L\times \R)$ with frame bounds $A \le B$, then for each $\delta > 0$ we have
 \begin{equation}\label{trace bounds}
A \ {\rm Tr\ }C \le T_\delta \le B \ {\rm Tr\ }C.
\end{equation}
and by Lemma \ref{trace class},
$$
{\rm Tr\ }C  =  \|h\|^2 \int_1^2 \int_N c(x,a) a^{-3} \ dxda = \|h\|^2 \ln 2.
$$
Hence it will suffice to show that as $\delta \rightarrow 0$,
\begin{equation}\label{suffice}
T_\delta \longrightarrow \frac{1}{\a \b}\  \|h\|^2 \ \int_{0}^\infty\int_\R \ |g(\l,t)|^2  \l^{-1}dt d\l.
\end{equation}
We begin by examining the expression $S(x,a)$ more closely.
Following \cite[page 65, proof of Theorem 3.3.1]{D}, for each $x_2, x_3, $ and $k$, we have
\begin{equation} \label{x3 ident}
\sum_{m \in \Bbb Z} \ w\left(\frac{|x_3+m+kx_2|}{a}\right) = a + \rho_3(x_3+kx_2,a)
\end{equation}
where $|\rho_3(x_3+kx_2,a)| \le \delta$ holds for all $k,x_2,x_3,$ and $a$. Similarly
\begin{equation} \label{x2 ident}
\sum_{l \in \beta\Bbb Z} \ w\left(\frac{|x_2+l|}{a}\right) = \frac{a}{\beta} + \rho_2(x_2,a)
\end{equation}
and
\begin{equation} \label{x1 ident}
\sum_{k \in \beta\Bbb Z} \ w\left(\frac{|x_1+k|}{a}\right) = \frac{a}{\a} + \rho_1(x_1,a)
\end{equation}
where $| \rho_2(x_2,a) |\le \delta$ and $ |\rho_1(x_1,a) |\le \delta$ hold for all $x_1, x_2$ and $a$. Applying (\ref{x3 ident})
 to the expression $S_\delta(x,a)$, we have
$$
\begin{aligned}
S_\delta(x,a) &=  \sum_{k\in\a\Bbb Z, l\in\b\Bbb Z} \  w\left(\frac{|x_1+k|}{\sqrt{a}}\right)w\left(\frac{|x_2+l|}{\sqrt{a}}\right) \sum_{m}w\left(\frac{|x_3+m+kx_2|}{a}\right)\\
&=  \sum_{k,l} \  w\left(\frac{|x_1+k|}{\sqrt{a}}\right)w\left(\frac{|x_2+l|}{\sqrt{a}}\right)\bigl(a + \rho_3(x_3+ky,a)\bigr)\\
&= a \ \sum_{k,l} \ w\left(\frac{|x_1+k|}{\sqrt{a}}\right)w\left(\frac{|x_2+l|}{\sqrt{a}}\right) + R_3(x,a)
\end{aligned}
$$
where
$$
R_3(x,a) = \sum_{k,l} \ w\left(\frac{|x_1+k|}{\sqrt{a}}\right)w\left(\frac{|x_2+l|}{\sqrt{a}}\right) \rho_3(x_3+kx_2,a).
$$
Now applying (\ref{x2 ident}) and (\ref{x1 ident}), we get
$$
S_\delta(x,a) = \frac{a^2}{\a \b} + R_\delta(x,a)
$$
where $R_\delta = R_{12} + R_3$ and
$$
R_{12}(x,a) = \frac{a \sqrt{a}}{\b} \rho_1(x_1,\sqrt{a}) +  \frac{a \sqrt{a}}{\a} \rho_2(x_2,\sqrt{a}) + a \rho_1(x_1,\sqrt{a})\rho_2(x_2,\sqrt{a}).
$$
So now we can write
$$
T_\delta =  \frac{1}{\a \b} \ \int_{0}^{\infty} \int_N \ \left|\langle g,  \hat T_{x} \hat D_{a} h\rangle \right|^2 \  \left(\frac{a^2}{\a \b} + R_\delta(x,a)\right) a^{-3} \ dx da
$$
We claim that
\begin{equation}\label{key int}
\int_{0}^{\infty} \int_N \ \left|\langle g,  \hat T_{x} \hat D_{a} h\rangle \right|^2 \ a^{-1}dx da < \infty.
\end{equation}
Observe that
$$
\begin{aligned}
|R_3(x,a)| &\le \sum_{k,l} \ w\left(\frac{|x_1+k|}{\sqrt{a}}\right)w\left(\frac{|x_2+l|}{\sqrt{a}}\right) |\rho_3(x_3+kx_2,a) | \\
&\le \delta \ \sum_{k,l} \ w\left(\frac{|x_1+k|}{\sqrt{a}}\right)w\left(\frac{|x_2+l|}{\sqrt{a}}\right)\\
&\le \delta \left(  \frac{\sqrt{a}}{\b} \rho_1(x_1,\sqrt{a}) +  \frac{ \sqrt{a}}{\a} \rho_2(x_2,\sqrt{a}) +  \rho_1(x_1,\sqrt{a})\rho_2(x_2,\sqrt{a})\right) \\
&\le \delta^2 \left(\frac{\sqrt{a}}{\b} + \frac{ \sqrt{a}}{\a} \right) + \delta^3
\end{aligned}
$$
while
$$
|R_{12}(x,a)| \le \delta \left(\frac{a\sqrt{a}}{\b} + \frac{a \sqrt{a}}{\a} \right) + \delta^2 a
$$
and so
$$
|R_\delta(x,a)| \le \delta \left(\frac{a\sqrt{a}}{\b} + \frac{a \sqrt{a}}{\a} \right) + \delta^2 \left(a + \frac{\sqrt{a}}{\b} + \frac{ \sqrt{a}}{\a} \right) + \delta^3.
$$
Hence $a^{-2} |R_\delta(x,a)| \rightarrow 0$ as $a \rightarrow \infty$. It follows that there is $m > 0$ such that for $a > m$,
$$
S_\delta(x,a) \ge \frac{a^2}{2\a\b}
$$
and hence
$$
\int_{m}^{\infty} \int_N \ \left|\langle g,  \hat T_{x} \hat D_{a} h\rangle \right|^2 \ a^{-1}dx da \le 2 \a\b \int_{m}^{\infty} \int_N \
\left|\langle g,  \hat T_{x} \hat D_{a} h\rangle \right|^2 \ S(x,a) \ a^{-3}dx da < \infty.
$$
Since $(x,a) \mapsto \left|\langle g,  \hat T_{x} \hat D_{a} h\rangle \right|^2 \ a^{-3}$ is integrable then the integral
 (\ref{key int}) taken from $a = 0$ to $m$ is  finite, and the claim follows. Now we can write
$$
T_\delta =  \frac{1}{\a \b} \ \int_{0}^{\infty} \int_N \ \left|\langle g,  \hat T_{x} \hat D_{a} h\rangle \right|^2 \  a^{-1} \ dx da  + I_\delta
$$
where
$$
I_\delta =  \int_{0}^{\infty} \int_N \ \left|\langle g, \hat T_{x} \hat D_{a} h\rangle \right|^2 \ R_\delta(x,a)  \ a^{-3} \ dx da.
$$
Since $ (x,a) \mapsto  \left|\langle g, \hat T_{x} \hat D_{a} h\rangle \right|^2$ is integrable over $(0,\infty)\times N$
 with respect to both of the measures $a^{-1} dx da$ and $a^{-3} dx da$, then it is integrable with respect to the measure
  $a^{-p} dx da$ for any $1 < p < 3$. So
$$
\begin{aligned}
I_\delta &=  \int_{0}^{\infty} \int_N \ \left|\langle g, \hat T_{x} \hat D_{a} h\rangle \right|^2 \ R_\delta(x,a)  \ a^{-3} \ dx da \\
&\le  \int_{0}^{\infty} \int_N \ \left|\langle g, \hat T_{x} \hat D_{a} h\rangle \right|^2 \ \left(\delta \left(\frac{a\sqrt{a}}{\b} + \frac{a \sqrt{a}}{\a} \right) + \delta^2 \left(a + \frac{\sqrt{a}}{\b} + \frac{ \sqrt{a}}{\a} \right) + \delta^3\right)  \ a^{-3} \ dx da \\
&= \delta\int_{0}^{\infty} \int_N \ \left|\langle g, \hat T_{x} \hat D_{a} h\rangle \right|^2 \  \left(\frac{1}{\b} + \frac{1}{\a} \right)  ( a^{-3/2} + \delta a^{-5/2}) \ dx da \\
&\hspace{1in} + \delta^2 \int_{0}^{\infty} \int_N \ \left|\langle g, \hat T_{x} \hat D_{a} h\rangle \right|^2 \   \ a^{-2} \ dx da \\
&\hspace{2in}+ \delta^3 \int_{0}^{\infty} \int_N \ \left|\langle g, \hat T_{x} \hat D_{a} h\rangle \right|^2 \ a^{-3} dx da
\end{aligned}
$$
and now we have
$$
T_\delta \longrightarrow \frac{1}{\a \b} \
\int_{0}^{\infty} \int_N \ \left|\langle g,  \hat T_{x} \hat D_{a} h\rangle \right|^2 \  a^{-1} \ dx da
$$
as $\delta \rightarrow 0$.

Finally,  let $\phi$ and $\psi$ be the preimages in $\mathcal H$ of $g$ and $h$.
 A standard calculation (see for example  \cite[equation (2.1)]{C3}) shows
$$
\begin{aligned}
 \int_{0}^{\infty} \int_N \ \left|\langle g, \hat T_{x} \hat D_{a} h\rangle \right|^2 \  a^{-1} \ dx da &=  \int_{0}^{\infty} \int_N \ \left| \phi * (\tau(a)\psi)^*(x) \right|^2 \ dx  a^{-1} da \\
 &= \int_{0}^{\infty}  \int_\L \ \| \hat \phi(\l) C_{a}^{-1} \hat \psi (a \l)^* C_a\|^2 |\l|d\l   a da.
 \end{aligned}
 $$
 where the norm in the last integral is the Hilbert-Schmidt norm. Now $\hat\phi(\l) = g(\l,\cdot) \otimes e_\l$, $\hat\psi(\l)^* = e_\l \otimes h(\l,\cdot)$. It follows that for a.e. $\l$,
\begin{equation}\label{ft expression}
  \| \hat \phi(\l) C_a^{-1} \hat \psi (a \l)^* C_a^{-1}\|^2 = \|g(\l,\cdot)\|_{L^2(\R)} \|h(a\l,\cdot)\|_{L^2(\R)}.
\end{equation}
Recalling that $h$ is supported on $(0,\infty) \times \R$, we see that (\ref{ft expression}) vanishes when $\l < 0$ and so
 $$
\begin{aligned}
 \int_{0}^{\infty} \int_N \ \left|\langle g, \hat T_{x} \hat D_{a} h\rangle \right|^2 \  a^{-1} \ dx da
 &= \int_{ 0}^\infty\int_{0}^\infty \|g(\l,\cdot)\|^2_{L^2(\R)} \ \|h(a \l,\cdot)\|^2_{L^2(\R)} \  \l d\l \ a da \\
  &=  \int_{ 0}^\infty\int_{0}^\infty\|g(\l,\cdot)\|^2_{L^2(\R)} \ \|h(a,\cdot)\|^2_{L^2(\R)} \ a da  \ \l^{-1}d\l \\
 &= \|h\|^2 \ \int_{0}^\infty\int_\R \ |g(\l,t)|^2  \l^{-1} dt d\l.
 \end{aligned}
 $$
 Thus the limit (\ref{suffice}) is shown, and the proof is complete.

\end{proof}

The following is immediate in light of \cite[Proposition 2.6]{C3}.

\begin{cor} Let $\psi$ belong to an invariant multiplicity one subspace $\mathcal H$. If $\psi$ is a Parseval frame vector for the unitary system $\mathcal U(\a,\b,\mathcal H) $, then $\frac{1}{\sqrt{\a\b \ln 2}}\psi$ is an admissible vector for the corresponding subrepresentation $\tau_\mathcal H$ of $\tau$.
\end{cor}

 Let $V : \mathcal H \rightarrow L^2(\L\times \R)$ be an invariant multiplicity-one subspace of $L^2(N)$, and let $I$ be a subset of $\L$. We will say that a function $g\in L^2(I\times \R)$ is an $I$-wavelet if $\hat{\mathcal U}(g,1,1,1)$ is a Parseval frame for $L^2(\L\times\R)$; we shall also say that $\psi = V^{-1}g$ is an $I$-wavelet.

\begin{lemma} \label{1 enough} Suppose that every Gabor field $g$ over $I$ is an $I$-wavelet.  Then we have the following.

\vspace{.1in} 
\noindent 
(a) $I$ is dilation congruent with the Shannon set.

\vspace{.1in}
\noindent
(b) For any lattice parameters $\a$ and $\b$ with $I \subset [-1/\a\b,1/\a\b]$, and for any unitary multiplier $\chi$,  $\widehat{\mathcal U}(g, \a,\b,\chi)$ is a Parseval frame for $L^2(\L\times \R)$.

\end{lemma}

\begin{proof} Let $g$ be a Gabor field over $I$. To show (a), we first observe that  $\R = \cup_{j\in\Bbb Z} 2^jI$. Otherwise, we would have a measurable $J \subset \R$ of finite (non-zero) measure such that $J \cap \bigl(\cup_{j\in\Bbb Z}2^jI\bigr) = \emptyset$. Then for all integers $j, k, l$ and $m$,
  $$
  \int_{\R} \int_{\L} \ \bold 1_{J\times [0,1]}(\l,t) \hat D_j \hat T_{k,l,m}g(\l,t) \ |\l|d\l dt = 0
  $$
contradicting our assumption about $g$. Now just as in \cite[proof of Theorem 5.4]{HL}, we have a measurable subset $K$ of $ I$ such that $K$ is dilation congruent with the Shannon set, so that $\R = \cup_{j\in\Bbb Z} 2^j K$ as a disjoint union, and we claim that $K = I$. Put $h := \bold 1_{K\times [0,1]}$; since  $K$ is dilation congruent with the Shannon set, then
$$
\int_0^\infty \int_\R \ |h(\l,t)|^2 \l^{-1}d\l dt =  \int_{K\cap (0,\infty)}\l^{-1} d\l
= \int_{[1/2,1]}\l^{-1} d\l = \ln 2.
$$
Since  $\widehat{\mathcal U}(g,1,1,1)$ is a Parseval frame for $L^2(\L\times \R)$,  then by Theorem \ref{nec condition}, we have
$$
\int_0^\infty \int_\R \ |g(\l,t)|^2 \l^{-1}d\l dt = \ln 2 .
$$
Since $g$ is a Gabor field we have $\int |g(\l,t)|^2 dt =1$ for a.e. $\l\in I$  from   (\ref{bandlimit}), and hence
$$
 \ln 2  = \int_0^\infty \int_\R \ |g(\l,t)|^2 \l^{-1}d\l dt = \int_{I\cap (0,\infty)}\l^{-1} d\l
$$

Similarly
$$
\int_{I\cap (-\infty,0)}\l^{-1} d\l 
= \ln 2 =  \int_{K\cap (-\infty,0)}\l^{-1} d\l
$$
It follows that $I = K$. Next we claim that $\widehat{\mathcal T}(g, 1,1)$ is a Parseval frame for $L^2(I\times \R)$. Let $f\in L^2(I\times \R)$ and recall that $\hat D_j f(\l,t) = f(2^j\l,2^{-j/2}t) 2^{3j/4}$, so that the support of $\hat D_jf$ is contained in $2^{-j}I \times \R$. Hence for $j \ne 0$,
$$
\langle f, \hat D_j \hat T_{k,l,m}g\rangle = \langle\hat D_{-j} f, g\rangle = 0, $$
and we have
$$
\|f\|^2 = \sum_{j,k,l,m} |\langle f, \hat D_j \hat T_{k,l,m}g\rangle|^2 =  \sum_{k,l,m} |\langle f, \hat T_{k,l,m}g\rangle|^2.
$$
Now since $g$ was an arbitrary Gabor field, we have that $I$ is a GH set. Part (b) now follows from  Propositions \ref{GHS} and \ref{the concrete frame}. 

\end{proof}

 The preceding result suggests a definition for wavelet sets in our setting that is analogous to the definition of wavelet sets (or frame sets) for Euclidean spaces.

\begin{definition} Let $I$ be a measurable subset of $[-1,1]$. We shall say that $I$ is a Heisenberg wavelet set if every Gabor field over $I$ is an $I$-wavelet.

  \end{definition}

The characterization of Heisenberg wavelet sets is exactly analogous with the classical case.

\begin{thm} \label{PW sets} Let $I$ be a measurable subset of $[-1,1]$.  Then the following are equivalent.

\vspace{.1in}
\noindent
(i) $I$ is a Heisenberg wavelet set.

\vspace{.1in} \noindent (ii) $I$ is translation congruent with a subset  of the unit interval and $I$ is dilation congruent with the Shannon set.

\end{thm}

\begin{proof}  If $I$ is a Heisenberg wavelet set then (ii) follows from Lemma \ref{1 enough} and Proposition \ref{GHS}. Suppose that (ii) holds and let $g$ be any Gabor field over $I$. By Proposition \ref{PF on I times R}, $\widehat{\mathcal T}(g, 1, 1)$ is a Parseval frame for $L^2(I\times \R)$. Now by Proposition \ref{the concrete frame}, $\widehat{\mathcal U}(g, 1, 1,1)$ is a Parseval frame for $L^2(\L\times \R)$.

\end{proof}

\vspace{1cm}

 Bradley Currey,
 {Department of Mathematics and Computer
Science, Saint Louis University, St. Louis, MO 63103}\\
 { \footnotesize{E-mail address: \texttt{{ curreybn@slu.edu}}}\\}

  Azita Mayeli,
  {Mathematics Department, Stony Brook University,  Stony Brook NY, 11794-3651, USA }\\
   \footnotesize{E-mail address: \texttt{{amayeli@math.sunysb.de}}}\\

\end{document}